\documentstyle[12pt]{article}
\newcommand{\const}{\mathop{\rm const}\limits}

\newcommand{\supp}{\mathop{\rm supp}\limits}

\begin{document}

\begin{center}

{\bf A BANACH REARRANGEMENT NORM CHARACTERIZATION  }\\

\vspace{3mm}

{\bf FOR TAIL BEHAVIOR  OF  MEASURABLE FUNCTIONS} \\

\vspace{3mm}

{\bf (RANDOM VARIABLES).  }

\vspace{4mm}

 $ {\bf E.Ostrovsky^a, \ \ L.Sirota^b } $ \\

\vspace{4mm}

$ ^a $ Corresponding Author. Department of Mathematics and computer science, Bar-Ilan University, 84105, Ramat Gan, Israel.\\
\end{center}
E - mail: \ galo@list.ru \  eugostrovsky@list.ru\\
\begin{center}
$ ^b $  Department of Mathematics and computer science. Bar-Ilan University,
84105, Ramat Gan, Israel.\\
E - mail: \ sirota@zahav.net.il\\

\vspace{3mm}
                    {\sc Abstract.}\\

 \end{center}

 \vspace{3mm}

 We construct a Banach rearrangement invariant norm on the measurable space
for which the finiteness of this norm for  measurable function (random variable) is equivalent to
suitable tail  (heavy tail and light tail) behavior.\par
 We investigate also a conjugate to offered spaces and obtain some embedding theorems.\par

 Possible applications: Functional Analysis (for instance, interpolation of operators), Integral Equations,
Probability Theory and Statistics (tail estimations for random variables) etc.\par

 \vspace{3mm}

 \par

\vspace{3mm}


\vspace{3mm}

\section{Notations. Statement of problem.}

\vspace{3mm}

 Let $ (X = \{x\}, \cal{A}, \mu) $ be measurable space with  non-trivial sigma-finite measure $ \mu. $
We will suppose without loss of generality in the case  $ \mu(X) < \infty $  that $ \mu(X) = 1 $ (the probabilistic case)
and denote $  x = \omega, \ {\bf P} = \mu. $\par
 Define as usually for arbitrary measurable function $ f: X \to R $ its distribution function
 (more exactly, tail function)

 $$
 T_f(t) = \mu\{x: |f(x)| \ge t \}, \  t \ge 0,
 $$

 $$
 ||f||_p = \left[ \int_X |f(x)|^p \ \mu(dx) \right]^{1/p}, \ p \ge 1; \ L_p = \{f, ||f||_p < \infty\},
 $$
 and denote by $ f^*(t) = T_f^{-1}(t) $ the left inverse to the tail function $ T_f(t);  $
$$
 f^{**}(t) \stackrel{def}{=} t^{-1} \int_0^t f^*(s) \ ds, \ t > 0.
$$

  We will denote the set of all tail functions  as $ \{ T \}; $ obviously, the set  $ \{ T \}$ contains on
all the functions $ \{ H = H(t), \ t \ge 0 \}  $ which are right continuous, monotonically non-increasing
with values in the set  $ [0, \mu(X)]. $ \par

 Let  $ w = w(s), s \ge 0 $ be any continuous strictly increasing   numerical function (weight) defined on the set
 $ s \in (0, \mu(X)) $  such that
 $$
 w(s) = 0 \Leftrightarrow s = 0;  \ \lim_{s \to \mu(X)}w(s) = \infty.\eqno(1.1)
 $$

 {\it In what follows the variables $ s, T = T(t) $ changes in the interval $ 0 < s,T < \mu(X).  $  \par
  Moreover, we impose on the set of all such a functions $ W = \{ w \} $ the following  restriction:  }\par

$$
\forall w \in W \ \exists T \in \{ T \}  \Rightarrow w(T(s)) = 1/s. \eqno(1.2)
$$

 Let us introduce the following important functional

 $$
 \gamma(w) = \sup_{t > 0} \left[  \frac{w(t)}{t} \ \int_0^t \frac{du}{w(u)} \right] \eqno(1.3)
 $$
and the following quasi-norms:

$$
||f||_w^* = \sup_{t > 0} [w(t) \ f^*(t)], \eqno(1.4)
$$

$$
||f||_w = \sup_{t > 0} [w(t) \ f^{**}(t)], \eqno(1.5)
$$
 The necessary and sufficient condition for finiteness of the functional $ \gamma(w) $ see, e.g.
in the article \cite{Astashkin1}. \par

\vspace{3mm}
{\bf Remark 1.1.} Note that
\vspace{3mm}
 $$
 ||f||^*_w = \sup_{t > 0} [t w(T_f(t))], \eqno(1.6)
 $$
so that if $ ||f||^*_w \in (0,\infty),  $ then
$$
T_f(t) \le w^{-1}(||f||^*_w/t).
$$

 Therefore the functional $ f \to ||f||^* $ may called "the tail quasinorm". \par
 Analogous  functionals was introduced in the books
\cite{Grafakos1}, chapter 1; \cite{Lieb1}, chapter 9;
in the  articles \cite{Liu1}, \cite{Liu2}, \cite{Liu3} etc.
 For instance,  in the  article \cite{Liu1} was introduced and investigated the functionals

 $$
 \rho_{\Phi}(f) = \sup_{t > 0} [ \Phi(t) \ T_f(t)  ], \eqno(1.7)
 $$

$$
||f||_{wL\Phi} = \inf \{ c, c > 0, \sup_{t > 0} [ \Phi(t/c) T_f(t) ] \le 1  \} \eqno(1.8)
$$
and the spaces $ wL\Phi= \{ f: ||f||_{wL\Phi} < \infty  \}, $

$$
wM\Phi = \{ f:  \forall c > 0 \  \Rightarrow \sup_{t > 0}  [ \Phi(t/c) T_f(t)  ] < \infty  \}. \eqno(1.9)
$$
 Here $ \Phi(\cdot)  $ is arbitrary Young-Orlicz function. \par

 The spaces $ wM\Phi, wL\Phi  $ was named in  \cite{Cwikel1}, \cite{Cwikel2}, \cite{Liu1} as
"weak  Orlicz spaces," in \cite{Grafakos1}, \cite{Lieb1} as "weak Lebesgue spaces".
 The  functionals of a type (1.5) - (1.9) are called "Lorentz weak norm" or "Marcinkiewicz norm",
see \cite{Astashkin1}, \cite{Carro1}, \cite{Krein1}, \cite{Lorentz1}, \cite{Soria1} etc. \par

 In the article \cite{Carro1}  was considered more general case of a functionals of a view (in our notations)

 $$
 ||f||_{w,p} = \left[ \int_9^{\infty} f^*(t) \ w(t) \ dt \right]^{1/p},
 $$
 $ ||f||_{w,\infty} = \sup_t [f^*(t) w(t)] $ and was obtained in particular the condition for quasi-normalizing of this spaces.\par

  For instance, if $ \Phi(t) = \Phi_p(t) = t^{1/p}, p \in [1,\infty),  $ then both the spaces
$ wM\Phi, wL\Phi  $ coincides with the Lorentz space $ L(p,\infty), $  see \cite{Bennet1}, chapter 4,
section 4, p. 216-217. \par
\vspace{3mm}
{\bf Remark 1.2.} As long as
\vspace{3mm}
$$
f^{**}(t)= t^{-1} \sup_{\mu(E) \le t} \int_E |f(x)| \ \mu(dx),
$$
we can rewrite  the expression for $ ||f||_w  $ as follows:

$$
||f||_w = \sup_{t > 0} \left[ (w(t)/t) \cdot \sup_{E: \mu(E) \le t} \int_E  |f(x)| \ \mu(dx)   \right]. \eqno(1.10)
$$

 If the measure $ \mu $ has not atoms, then the expression (1.10) may be rewritten as follows:

 $$
||f||_w = \sup_{ E: 0 < \mu(E) < \infty } \left[ \frac{ w(\mu(E))}{\mu(E)} \cdot \int_E |f(x)| \ \mu(dx) \right]. \eqno(1.10a)
 $$

  It follows from  equality (1.10) that $ ||f||_w $ is true rearrangement invariant norm and the space
$ L_w = \{ f: ||f||_w < \infty  \}  $ is complete Banach functional rearrangement invariant
 space with Fatou property.  The proof is alike to one in the case $  w(t) = t^{1/p}, \ p \ge 1;  $ see \cite{Bennet1}, chapters 1,2;
\cite{Stein1}, chapter 8. \par

\vspace{3mm}

 The norm $  ||f||_w   $ is named Marcinkiewicz's norm, see \cite{Krein1},  chapter 2, section 2. \par

\vspace{3mm}

{\bf Example 1.1. Fundamental function.}\par
 Let  $ \delta  $ be arbitrary number from the set $ (0, \mu(X) ) $ and let $ B $ be any measurable set such that
  $ \mu(B) = \delta. $  The function $ \phi(\delta; Z) $ defined aside from $  \delta $ on the rearrangement
 invariant space $  Z $  as follows:

  $$
\phi(\delta; Z)  = ||I(B)||Z
  $$
is called a fundamental function of the space $ Z. $ \par

 If the measure $ \mu $ has not atoms, then the formula  (1.10a) gives us:

 $$
\phi(\delta; L_w) = w(\delta).
 $$

\vspace{3mm}
{\bf Remark 1.3.}  The equality (1.6)  may be used for the following definition. Let $ g $ be absolutely integrable
non-zero function: $ g \in L_1(X, \mu). $ We define the {\it natural} weight function $ w^{(g)}= w^{(g)}(s) $ as follows:
\vspace{3mm}
 $$
 1= t \cdot w^{(g)}(T_g(t)), \eqno(1.11)
 $$
or equally
$$
w^{(g)}(s) = T_g^{-1}(1/s). \eqno(1.12)
$$
 Note that $ ||g||^*_{ w^{(g)} } = 1.  $ \par
 We will say that in this case the function $ g(\cdot) $ generated  the correspondent space $  L_{w^{(g)}} $ and
 $  L_w^{(g,*)} = \{ f: ||f||^*_{w^{(g)} < \infty} \}. $ \par
\vspace{3mm}

{\bf  Our aim in this short report is to prove (under simple conditions) that the quasinorm $ ||f||^*_w  $ and the norm
$ ||f||_w  $ are linear equivalent. }\par
 This fact is true for heavy tails;  the case of light tails will be considered further, in which we prove that 
the weak Orlicz's norm is equivalent ordinary  Orlicz's norm for exponential correspondent Young function. \par 

\vspace{3mm}

  Our results improve ones ib \cite{Astashkin1},
 \cite{Krein1},  chapter 2, section 2;   see also reference therein. In
 particular, we find the exact value for estimated functional (Theorem 2.1), consider the case of light
 tails (Theorem 3.1), generalize the embedding theorem in \cite{Kufner1}, p. 167   etc. \par

\vspace{3mm}

\section{Main result: the case of heavy tails.}

\vspace{3mm}

{\bf  Theorem 2.1.} Let
$$
w \in W, \ \gamma(w) < \infty,  \eqno(2.1)
$$
then

$$
1 \cdot||f||^*_w  \le ||f||_w \le \gamma(w) \cdot ||f||^*_w, \eqno(2.2)
$$
 and both  the coefficients $ "1" $ and $ "\gamma(w)" $ in (2.2) are the best possible. \par
\vspace{3mm}
{\bf Proof} is at the same as for the spaces $ L(p, \infty), $ see \cite{Stein1}, chapter 8.\par
\vspace{3mm}
{\bf A. Inequalities.}\\
\vspace{3mm}
 The left-hand side of assertion (2.2) follows immediately from the inequality $  f^{**}(t) \ge f^*(t) $
 even without the conditions (2.1). \par
  Let now  $ \gamma(w) < \infty.$ We have:

 $$
w(t) f^{**}(t) = w(t) \ \frac{1}{t} \ \int_0^t f^*(u) du = \frac{w(t)}{t} \ \int_0^t \frac{w(u) f(u) du}{w(u)} \le
 $$

$$
 \frac{w(t)}{t} \ ||f||^*_w \ \int_0^t \frac{du}{w(u)} \le \gamma(w) \ ||f||^*_w.\eqno(2.3)
$$

Taking supremum of the last inequality over all the values $ t, $ we obtain the right-hand side (2.2).\par
\vspace{3mm}

{\bf B. Exactness.} \\

\vspace{3mm}

 The exactness of the constant $ " 1 " $ follows immediately from the consideration of the case $ w = w_p, \ p > 1. $
Namely, in this case we obtain the classical inequality

$$
||f||^*_{w_p}  \le ||f||_{w_p} \le \frac{p}{p-1} \cdot ||f||^*_{w_p};
$$
note that  $ \lim_{p \to \infty} p/(p-1)=1.  $\par
 It remains to prove the exactness of the coefficient $ \gamma(w) $ in the right-hand side of inequality (2.2). It is reasonable to
suppose $ \gamma(w) < \infty; $ in other case it is nothing to prove. \par
 Let here $  w \in W; $  {\it we can prove  moreover that the right-hand inequality in (2.2) is exact still for arbitrary function} $ w. $  \par
 There exists a positive (measurable) integrable function  $ g: X \to R_+ $  for which

 $$
 w(t) = w^{(g)}(t) = 1/T_g^{-1}(t) = 1/g^*(t), \eqno(2.4)
 $$
then  $ ||g||^*_{ w^{(g)} } = 1 $  and both the spaces $ L_w, \ L^*_w $ are generated by means of the function $ g. $ \par

 Denote for any function $  w \in W $ the following functional

$$
G(w) := \sup_{f: ||f||^*_w = 1} \left[ \frac{||f||_w}{||f||^*_w} \right] = \sup_{f: ||f||_w = 1} ||f||_w. \eqno(2.5)
$$

 We have:

 $$
 G(w) \ge ||g||_w = \sup_t \left[ w(t) g^{**}(t) \right] = \sup_t \left[ w(t) \cdot \frac{1}{t} \int_0^t g^*(s) ds \right] =
 $$

$$
\sup_t \left[ w(t) \cdot \frac{1}{t} \int_0^t \frac{ ds}{w(s)} \right] = \gamma(w),
$$
Q.E.D. \par

\vspace{3mm}
{\bf Example 2.1.}  {\it We suppose in all considered in this article examples that $ \mu = {\bf P} $ or equally} $ \mu(X) = 1. $\par

\vspace{3mm}

  Let $  w(s)= w_{p,l}(s) = s^{1/p} \ l(s), \ 1 \le p < \infty, $ where $ l = l(s) $ is non-negative positive for positive values $ s $
continuous on the semi-open interval  $ 0 < s \le 1, $  slowly varying  function  as $ s \to 0+.  $  By definition, $ w_p(s) = w_{p,1}(s) =
s^{1/p}. $ \par
    The condition (2.1) is satisfied iff $ p > 1. $\par
 The asymptotical as $ t \to 0+ $  relation

 $$
 I := \int_0^t \frac{du}{w_{p,l}(u)} \asymp C \frac{t}{w_{p,l}(t)}, \  p > 1 \eqno(2.6)
 $$
may be obtained from the book \cite{Bingham1}, chapter 1, sections 1.5., 1.6 p. 26-27. \par

 As a consequence: for a random variable $ \xi = \xi(\omega) $ the tail inequality

 $$
 T_{\xi}(t) \le C_1 K^p t^{-p} \log^{-\kappa p  }(t/K), \ t \ge 2 K,
 $$

 $$
 C_1,K = \const > 0, \ \kappa = \const, p = \const > 1 \eqno(2.7)
 $$
 is equivalent to the following norm estimation:

 $$
 \sup_{0 < t \le 1/2} \left\{t^{1/p-1} \ |\log^{\kappa}(t)| \ \sup_{ 0 < {\bf P}(E) \le t  }\int_E |\xi(\omega)| \ {\bf P}(d \omega) \right\} \le C_2(p,\kappa) \ K.
 \eqno(2.8)
 $$

 If the measure $ {\bf P} $ is atomless, then the inequality (2.8)  may be simplified as follows:

 $$
   \sup_{ 0 < {\bf P}(E) \le 1/2  } \left\{ [{\bf P}(E)]^{1/p-1} |\log^{\kappa} {\bf P}(E) | \
 \int_E |\xi(\omega)| \ {\bf P}(d \omega) \right\}  \le C_2(p,\kappa) \ K. \eqno(2.9)
 $$

 Notice that this norm description of heavy tail distributed random variables is more convenient as description by means of the
so-called  moment, or Grand Lebesgue spaces, as long as ones are not completely adequate, see the example 5.1 in the article
\cite{Ostrovsky1}.\par

\vspace{3mm}

\section{Main result: the case of light tails.}

\vspace{3mm}

{\bf Example 3.1.}

\vspace{3mm}

Let $  X = (0,1) $ with Lebesgue measure, and let

$$
h(x) = |\log \ x|, \ x \in X,  \eqno(3.1)
$$
then $ T_h(t) = e^{-t}, \ t \ge 0; $

$$
h^*(s) = |\log \ s|; \ w^{(h)}(s) = |\log \ s|^{-1}, \ s \in (0,1); \ ||h||^*_{w^{(h)}} = 1; \eqno(3.2)
$$
but

$$
f^{**}(t) = \frac{1}{t}\cdot \int_0^t |\log \ x| dx = 1 + |\log t|,
$$
and

$$
||h||_{w^{(h)}} = \sup_{t > 0} \left[w^{(h)}(t) \ h^{**}(t) \right] = \sup_{t \in (0,1) } \frac{1+|\log t|}{|\log t|} = \infty.\eqno(3.3)
$$

 Analogous implication $ ||h||^*_{w^{(h)}} = 1, \ ||h||_{w^{(h)}} = \infty $ if true for the functions of a view

 $$
 h(x) = h_{m,l(\cdot)}(x) = |\log x|^m \ l(|\log x|), \ m = \const > 0,
 $$
$ l(z) $ is slowly varying as $ z \to \infty $ non-negative function. \par
 Notice that the functions of a view $ h_{m,l(\cdot)}(x), \ x \in X $ have a light tails.  Indeed,  they belong to the
so-called {\it exponential} Orlicz spaces.\par

 Recall that the Young-Orlicz function $ N = N(u), \ u \in R $ is called {\it  exponential Young-Orlicz function,  }
 briefly: EOF, if

$$
N(u) = e^{\nu(u)} - 1, \eqno(3.4)
$$
where $ \nu =\nu(u) $ is even twice  continuous differentiable convex function, strictly monotonically  increasing on the right-hand
semi-axis and  such that
$$
 \nu(u) = 0 \ \Leftrightarrow u = 0; \ \nu'(0) = 0; \ \lim_{u \to \infty} \nu'(u) = \infty. \eqno(3.5)
$$

 The Orlicz space $  L(N) $  with  Young-Orlicz function $ N = N(u) $
 defined over Probabilistic space $ (X, \cal{A}, {\bf P}) $ is said to be
 {\it exponential Orlicz space,} briefly: EOS, if the function $  N = N(u) $ is exponential Orlicz function EOF.\par

 We denote in accordance with \cite{Kufner1}, p. 167 as $ F{N} $ the  set of all (measurable) functions $ u = u(x),
 u: X \to R   $ such that

 $$
 \forall k > 0 \ \Rightarrow \int_X N(k |u(x)|)  \ \mu(dx) < \infty. \eqno(3.6)
 $$

  Evidently, $ F(N) \subset L(N). $ \par

 Further, the subspace $ E(N)  $ of the whole space $ L(N) $
 is by definition the closure of the set of all bounded measurable functions
supported on the set of  finite measure.  It is known (see, e.g.  \cite{Kufner1}, p. 167)
that $ E(N) = F(N) . $\par

\vspace{3mm}

{\bf Theorem 3.1.}  Let the function $ \Phi = \Phi(u) $ be EOF,
so that the Orlicz space $ L(\Phi) $ is EOS. Then: \par
\vspace{3mm}
{\bf A.} The weak Orlicz space $ wL\Phi  $  coincides as the set equality with norm equivalence with
the ordinary (exponential) Orlicz space $ L(\Phi). $ \par
\vspace{3mm}
{\bf B.} The weak Orlicz space $ wM\Phi  $  coincides as the set equality with norm equivalence with
the subspace $  E(\Phi). $ \par

\vspace{3mm}

{\bf Proof \ A.} Let  $ f \in L(\Phi), \ f \ne 0; $  we can suppose without loss of generality $ ||f||L(\Phi) = 1. $
Then

$$
\int_X \Phi(|f(x)|) {\bf P}(dx) \le 1.
$$
 We use the Tchebychev's  inequality:

 $$
 {\bf P} ( |f| > C ) \le \frac{1}{\Phi(C)}, \ C = \const > 0;\eqno(3.7)
 $$
therefore $ f \in  wL\Phi $ and $ ||f||_{wL\Phi} \le 1. $ \par
 Conversely, let $ f \in  wL\Phi $ and $ ||f||_{wL\Phi} \le 1. $ It follows from the definition of the functional
 $  f \to ||f||_{wL\Phi} $  that

 $$
 T_f(t) \le \frac{1}{\Phi(t)}, \ t > 0
 $$
or equally

$$
T_f(t) \le C_1 e^{- \nu(t)  },  \ t \ge 1. \eqno(3.8)
$$
 It follows from the estimate  (3.8) that

$$
 \exists C_2 = C_2(C_1; \nu(\cdot)) \in (0,\infty) \ \Rightarrow
 \int_X  \Phi(|f(x)|/C_2) \ {\bf P}(dx) < \infty;\eqno(3.9)
$$
$$
    ||f||L(\Phi) \le C_3 (C_1, C_2) < \infty, \eqno(3.10)
$$
see \cite{Kozachenko1}; more detail explanation see in a monograph
\cite{Ostrovsky3}, chapter 1, section 1.2. \par
\vspace{3mm}
{\bf Proof \ B.} Let $ f \in wM\Phi; $  we conclude using at the same arguments as before

$$
 \forall C_4 > 0  \Rightarrow
 \int_X  \Phi(|f(x)|/C_4) \ {\bf P}(dx) < \infty. \eqno(3.11)
$$

 It follows from the  proposition (3.11) that the function $ f $ belongs to the subspace
$ F(\Phi); $  but we know $ F(\Phi) = E(\Phi). $ \par
 The converse inclusion $   E(\Phi) \subset wM\Phi $ follows from the simple verified fact that for every
measurable set $ B $  with finite measure $ \mu(B) < \infty $  its indicator $ I(B) $ belongs to the set
$ wM\Phi $ since

$$
T_{ I(B) }(t) = 0, \ t > 1.
$$

\vspace{3mm}

 Recall that the Absolutely Continuous Norm (ACN =  ACN(Z)) part of the rearrangement invariant (r.i.) space $ (Z, ||\cdot||Z) $
over the triple $ (X, \cal{A}, {\bf P}) $ consists on the functions $ \{ f \} $ with Absolutely Continuous Norm (ACN):

$$
\lim_{\epsilon \to 0+} \sup_{B \in \cal{A}, {\bf P}(B) < \epsilon } || f \cdot I(B)||Z  =0. \eqno(3.12)
$$
 The $ I(B) $ denotes the indicator function of the (measurable) set $ B. $
 The ACN part is also closed subspace of the space $ (Z, ||\cdot||Z) $. \par
More information about the spaces $ E(N) $ and $ ACN(Z) $ see in the classical monograph belonging to
C.Bennett and R.Sharpley \cite{Bennet1}, chapter 1, sections 2,3.  In particular, it is proved that
(in our notations) that $  ACN(L(\Phi)) \subset E\Phi. $ As a consequence: \par
\vspace{3mm}
 {\bf Remark 3.1.}  Every function from the set  $ wM\Phi $ has absolutely continuous norm in the whole space
 $ L(\Phi). $ \par
 As a contradiction: \par
\vspace{3mm}
 {\bf Remark 3.2.}  Both the spaces: $ wL\Phi,  $ and $ L_w  $ under conditions (2.1) have not ACN property. \par
  Indeed, the space $ wL\Phi  $ coincides with the Orlicz space $ L(\Phi)  $ with the Young-Orlicz function
 not satisfying the $  \Delta_2  $  condition. \par
  It remains to consider the space $ L_{w_p}(0,1), \ p > 1  $ with Lebesgue measure $ m. $
  Choose a function $ f(x) = x^{-1/p}  $ and let
 $ b \downarrow 0, \ b < 1; $  then $ m \{ (0,b) \} = b \downarrow 0; $

$$
||f \cdot I(0,b)||^*_w = \sup_{ t \in (0,1)} \left[ \frac{1}{t^{ 1-1/p}}\int_0^{\min(b,t) } x^{-1/p} \ dx \right] =
$$
$$
\frac{p}{p-1} \sup_{ t \in (0,1)} \left[ \frac{ (\min(b,t)^{1-1/p} )}{t^{1-1/p} } \right] = \frac{p}{p-1}
$$
and

$$
\lim_{b \downarrow 0} ||f \cdot I(0,b)||^*_w  = \frac{p}{p-1} > 0.\eqno(3.13)
$$

\vspace{3mm}
 {\bf Remark 3.3.}
  The complete description of the conjugate (= associate or dual) space to the considered here weak spaces:
weak  Orlicz space $L(\Phi)  $ etc.
see in the book \cite{Rao1}, chapter 11.  The dual spaces to the  Marcinkiewicz and Lorentz spaces, conditions of its
reflexivity  and separability, description of compact subsets,
for instance, weak   $  L(p), \ p \ge 1 $ spaces are described, e.g. in the book  \cite{Krein1}, chapter 2, section 3; in the
articles \cite{Cwikel2}, \cite{Cwikel3} etc. \par

\vspace{3mm}
 {\bf Remark 3.4.}

Let $ w \in W $  and let the measure $ \mu $ be resonant; then the fundamental function  of the associate space  $ (L_w)' $ has a view

 $$
 \phi(\delta; (L_w)') = \frac{\delta}{w(\delta)}. \eqno(3.14)
 $$

\vspace{3mm}
 {\bf Remark 3.5.} Obviously, the assertion of theorem 3.1 remains true if instead the weak Orlicz norm  stated the 
 functional (1.7) $ \rho_{\Phi}(f) = \sup_{t > 0} [ \Phi(t) \ T_f(t)  ]. $ \par 
 
 \vspace{4mm}

\section{Embedding theorem.}

\vspace{3mm}

 In the book \cite{Kufner1},  p. 105, theorem 2.18 is proved the following embedding theorem. Let
 $  p > 1, \epsilon = \const \in (0, p-1);  $  then

 $$
 L_p \subset L_{w,p} \subset L_{p-\epsilon}. \eqno(4.1)
 $$

 Hereafter in this section the symbol $ \subset $ denotes the linear continuous embedding.\par
{\it   We intend here to generalize the assertion (4.1) on the case of general $ L_w $ spaces
over probabilistic measure space.}
  We consider  the case when for $ p_0 = \const > 1, \Delta = \const \ge 0, $

 $$
 \Phi(u) \stackrel{def}{=} \Phi_{p_0, \Delta,S}(u) = |u|^{p_0} \ (\log |u|)^{-\Delta} \ S(\log |u|), \  |u| \ge e,   \eqno(4.2)
 $$

 $$
   \Phi_{p_0, \Delta,S}(u) = C(p_0,\Delta) \ u^2, \ |u| \le e;  \  C(p_0,\Delta) e^2 = e^{p_0 } S(1),
 $$
where $ S = S(z), \ z \ge 1 $ is continuous positive slowly varying as $ z \to \infty $  function.\par
 Note that the function $ \Phi_{p_0, \Delta}(u) $  satisfies the $\Delta_2  $ condition; the case of exponential
Orlicz space $ L(\Phi) $  is here trivial.\par
 In order to formulate our result, we need to recall some facts about the so-called Grand Lebesgue spaces (GLS).
     Recently, see       \cite{Fiorenza3},  \cite{Kozachenko1},\cite{Liflyand1},  \cite{Ostrovsky1}    etc.
     appears the so-called Grand Lebesgue Spaces $ GLS = G(\psi) =G\psi =
    G(\psi; A,B), \ A,B = \const, A \ge 1, A < B \le \infty, $ spaces consisting
    on all the measurable functions $ f: T \to R $ with finite norms

$$
     ||f||G(\psi) \stackrel{def}{=} \sup_{p \in (A,B)} \left[ |f|_p /\psi(p) \right].
     \eqno(4.3)
$$

      Here $ \psi(\cdot) $ is some continuous positive on the {\it open} interval
    $ (A,B) $ function such that

     $$
     \inf_{p \in (A,B)} \psi(p) > 0, \ \psi(p) = \infty, \ p \notin (A,B).
     $$
We will denote
$$
 \supp (\psi) \stackrel{def}{=} (A,B) = \{p: \psi(p) < \infty, \} \eqno(4.4)
$$

The set of all $ \psi $  functions with support $ \supp (\psi)= (A,B) $ will be
denoted by $ \Psi(A,B). $ \par
  This spaces are rearrangement invariant and  are used, for example, in the theory of probability,
  theory of Partial Differential Equations,  functional analysis, theory of Fourier series, theory of martingales,
  mathematical statistics, theory of approximation    etc.\par

 We will use the following important example (more exactly, the {\it  families
of examples}) of the $ \psi $ functions and correspondingly the GLS spaces.\par
Let us denote
$$
\psi(B,\beta;p) \stackrel{def}{=} (B - p)^{-\beta},\eqno(4.5)
$$
where $ \beta = \const \ge 0, 1 < B < \infty; p \in [1,B) $ so that
$$
\supp \psi(B, \beta;\cdot) = [1,B).
$$

As a particular case: let us introduce the following $  \psi $ function:

$$
\psi_{p_0, \Delta, S }(p)= (p_0 - p)^{-(1+\Delta)/p_0} \ S^{1/p_0} \left( \frac{p_0}{p_0 - p} \right), \eqno(4.6)
$$
$ \Delta = \const > -1, \ p_0 = \const > 1; \ 1 \le p < p_0 $
and the following spaces: $  L(\Phi_{p_0, \Delta,S}), wL\Phi_{p_0, \Delta,S} $ and $ G\psi_{p_0, \Delta, S }. $

\vspace{3mm}

{\bf Remark 4.1.} If we define the {\it degenerate } $ \psi_r(p), r = \const \ge 1 $ function as follows:
$$
\psi_r(p) = \infty, \ p \ne r; \psi_r(r) = 1
$$
and agree $ C/\infty = 0, C = \const > 0, $ then the $ G\psi_r(\cdot) $ space coincides
with the classical Lebesgue space $ L_r. $ \par

\vspace{3mm}

{\bf Remark 4.2.} Let $ \xi: X \to R $ be some (measurable) function from the set
$ L(p_1, p_2) = \cup_{p \in (p_1, p_2)} L(p), \ 1 \le p_1 < p_2 \le \infty. $ We can introduce the so-called
{\it natural} choice $ \psi_{\xi}(p) $  as as follows:

$$
\psi_{ \xi}(p) \stackrel{def}{=} ||\xi ||_p; \ p \in (p_1,p_2).
$$

 Evidently, in the case when $ \mu(X) < \infty, \  $ then by virtue of Lyapunov's inequality
 $ L(p_1, p_2) = L[1, p_2). $ \par

\vspace{3mm}

 Analogously,  let $ \xi = \{ \xi(y) \}, \ y \in Y $ be some {\it family} of measurable function uniformly
from the set $ L(p_1, p_2) = \cup_{p \in (p_1, p_2)} L(p), \ 1 \le p_1 < p_2 \le \infty. $ We can introduce the so-called
{\it natural} choice $ \psi_{ \{\xi \} }(p)$  as as follows:

$$
\psi_{ \{\xi \}}(p) \stackrel{def}{=} \sup_{y \in Y} ||\xi(y)||_p; \ p \in (p_1,p_2).
$$

\vspace{4mm}

{\bf Theorem 4.1.}

$$
 L(\Phi_{p_0, \Delta,S}) \subset wL\Phi_{p_0, \Delta,S} \subset G\psi_{p_0, \Delta, S }. \eqno(4.7)
$$

\vspace{3mm}

{\bf Proof. \ 1. Left hand-side inclusion. } \par

\vspace{3mm}

 Let $ f \in L(\Phi_{p_0, \Delta,S}), \ f \ne 0;  $ then  for some finite positive constant $  u \in (0, \infty) $

 $$
 \int_X \Phi_{p_0, \Delta,S}(|f(x)|/u) \ {\bf P}(dx) \le 1.
 $$
We use the Tchebychev's inequality:

$$
{\bf P}(|f(x)|/u > C) \le   \frac{1}{\Phi_{p_0, \Delta,S}(C/u)}, \  C \in (0, \infty),
$$
or equally

$$
{\bf P}(|f(x)|> C_1) \le   \frac{1}{\Phi_{p_0, \Delta,S}(C_1)}, \  C_1 \in (0, \infty), \eqno(4.8)
$$
 The last inequality (4.8) imply that $ f \in wL\Phi_{p_0, \Delta,S}. $ \par

\vspace{3mm}
   {\bf 2. Right hand-side inclusion.} \par
\vspace{3mm}
Let $ f \in wL\Phi_{p_0, \Delta,S}.  $ We can and will suppose without loss of generality that

$$
T_f(t) \le  t^{-p_0} \ \log^{\Delta}(t) \ S(\log t), \ t > 2. \eqno(4.9)
$$

 Let $ 1 \le p < p_0. $  Recall that in this section $ {\bf P} (X) = 1, $
 therefore $ T_f(t) \le 1. $   As long as

$$
\int_X |f(x)|^p \ {\bf P}(dx) = p  \int_0^{\infty} t^{p-1} \ T_f(t) \ dt,  \eqno(4.10)
$$
we obtain substituting the estimate (4.9) into (5.10), denoting $ \epsilon = p_0 - p $ and
tacking into account that $ \epsilon \to 0+: $

$$
\int_X |f(x)|^p \ {\bf P}(dx) \le p_0 \ \int_0^2 2^{p_0 - 1} dt + \int_2^{\infty}t^{p-1} \ T_f(t) \ dt;
$$

$$
\int_X |f(x)|^p \ {\bf P}(dx) - C_2 \le \int_1^{\infty} t^{p-p_0-1} \ \log^{\Delta}(t) \ S(\log t) \ dt  =
$$

$$
\int_0^{\infty} e^{-\epsilon y} \ y^{\Delta} \ S(y) \ dy = \epsilon^{-1 - \Delta}\int_0^{\infty} e^{-x} \ x^{\Delta} \
S(x/\epsilon) \ dx \sim
$$

$$
 \epsilon^{-1 - \Delta} \ \int_0^{\infty} e^{-x} \ x^{\Delta} \ S(1/\epsilon) dx =
 \epsilon^{-1 - \Delta} S(1/\epsilon) \Gamma(1 + \Delta).
$$
 Thus,  we have as $ p \to p_0-0 $

$$
||f||_p \le C \  (p_0-p)^{ - (1 + \Delta)/p_0 } \ S^{1/p_0}(1/(p_0-p)) \Gamma(1+\Delta) \asymp
\psi_{p_0, \Delta, S }(p). \eqno(4.11)
$$

 The detail grounding of passage to limit see in \cite{Liflyand1}, \cite{Ostrovsky1}.\par
This completes the proof of theorem 4.1. \par
\vspace{3mm}

 In order to show the exactness of proposition  the theorem 4.1, let us consider the following
inverse result. \par

\vspace{3mm}

{\bf Theorem 4.2.} We assert under at the same conditions as in theorem 4.1 that

$$
  G\psi_{p_0, \Delta, S } \subset wL\Phi_{p_0, \Delta + 1,S}. \eqno(4.12)
$$

\vspace{3mm}

{\bf Proof.} Let $ f \in   G\psi_{p_0, \Delta, S }; $ this imply by definition of the norm in
Grand Lebesgue Spaces for the values $ 1 \le p < p_0 $

$$
||f||_p \le  C \ \psi_{p_0, \Delta, S }(p)= C (p_0 - p)^{-(1+\Delta)/p_0} \ S^{1/p_0} \left( \frac{p_0}{p_0 - p} \right). \eqno(4.13)
$$
 We can  and will suppose in (4.13) without loss of generality  $ C = 1. $\par

  By means of Tchebychev's  inequality

 $$
 T_f(t) \le \frac{ \psi^p_{p_0, \Delta, S }(p)}{t^p}, \ t > 0,  \eqno(4.14)
 $$
therefore

 $$
 T_f(t) \le  \inf_{ 1 \le p < p_0 } \left[\frac{ \psi^p_{p_0, \Delta, S }(p)}{t^p} \right] \le
 $$

$$
 C_2 \  t^{-p_0} \ \log^{\Delta+1}(t) \ S(\log t), \ t > 2. \eqno(4.15)
$$

 The last inequality implies that $ f \in wL\Phi_{p_0, \Delta + 1,S},  $  Q.E.D. \par
 Notice that the estimate (4.15) is in general case non-improvable, see
\cite{Ostrovsky8}. \par

\vspace{3mm}
{\bf Remark 4.3.}  The  right-hand side inequality (4.7) of theorem 4.1 is non-improvable. 
 Namely, for the random variable  $ \eta, $  defined on the sufficiently rich probabilistic space
with the tail behavior 

$$
T_{\eta}(t) \sim  t^{-p_0} \ \log^{\Delta}(t) \ S(\log t), \ t \to \infty \eqno(4.16)
$$
we have  as $  p \to p_0 -0 \  \Rightarrow  ||\eta||_p \asymp $

$$
\Gamma(1+\Delta) \ (p_0-p)^{ - (1 + \Delta)/p_0 } \ S^{1/p_0}(1/(p_0-p)) = \Gamma(1+\Delta) \ 
\psi_{p_0, \Delta, S }(p). \eqno(4.17)
$$

\end{document}